# Yet another proof of Tychonoff's Theorem[1]

N. Noble[2]

**Abstract**

In 1971 I announced a "nice" proof of Tychonoff's Theorem, an immediate corollary of a result concerning closed projections combined with Mrówka's characterization of compactness: a space X is compact if and only if for each space Y the projection $\pi_Y: X \times Y \to Y$ is closed. I described the proof as "to appear" but to date it has not. In 2019 I published a generalization of the stronger closed projection result which yields a different look to the "nice" proof. Both versions are presented here.

## 0  Introduction

In Noble [1971] I suggested that "the proof of Tychonoff's theorem using Theorem 1.8 is also of some interest. Its use of the axiom of choice differs from that of other proofs in that it is restricted to statements about the index set" and referred the reader to "A nice proof of Tychonoff's Theorem (to appear)". That proof did not "appear"; the draft article was never finished and no longer exists. My vague recollection of it is that it used the proof as a springboard for examining variations of the Axiom of Choice. All that remains, aside from that recollection, is a footnote in Halpern, [1972, p119] which mentions, with regard to its relationship to some other conditions of the statement "Every infinite set algebra has a non-principal prime ideal.", that "The referee has informed us that this proof was discovered independently by Norman Noble."[3].

Fast forward: recently Matheron[2020] published three proofs of Tychonoff's Theorem, one of which relies upon Mrówka's characterization of compactness (and triggered my recollection of the subject of this note). Matheron's proof was extracted from a more general categorical result in Clementino & Tholen [1996][4], as was a similar proof given by Todd Trimble in Mathoverflow 26416 [2016]. Another proof using Mrówka's characterization is contained in Jachymski [2005], where some of the other consequences of my Theorem 1.8 are rediscovered.

## 1  The "nice" proof

For $\pi_Y$ the projection from $X \times Y$ onto $Y$ my 1971 proof exploits Mrówka's[5] characterization of compactness in the form: for S a subset of $X \times Y$, if $\pi_Y S$ is not closed then S is not closed. The proof is a transfinite induction which, for $\pi_Y S$ not closed, constructs a point x in clS\S. The induction step is a barely deeper characterization of closed projections, stated next.

**1.1 Proposition.**  The projection $\pi_Y: X \times Y \to Y$ is closed if and only if for each subset S of $X \times Y$ and each point y of $cl(\pi_Y S)$ there exists a point x of Y such that (x,y) is in cl(S).

---

**Proof.** First suppose $\pi_Y$ is closed, $S \subseteq X \times Y$, and $y \in \text{cl}(\pi_Y S)$. Since $\pi_Y \text{cl}(S)$ is a closed set containing $\pi_Y S$, $y$ is in $\pi_Y \text{cl}(S)$, that is, $X \times \{y\}$ meets $\text{cl}(S)$. For any point $x$ in that intersection, $(x,y)$ is in $\text{cl}(S)$. For the converse, if $\pi_Y$ is not close there exists a closed subset $S$ of $X \times Y$ with a point $y$ in $\text{cl}(\pi_Y S) \setminus \pi_Y S$ and for that $y$, no point $(x,y)$ is in $\text{cl}(S) = S$.

**Notation.** For $X = \Pi_{\alpha \in \kappa} X_\alpha$, $\kappa$ a cardinal, and $\beta \leq \kappa$, set $X^\beta = \Pi_{\alpha < \beta} X_\alpha$, so $X^{\beta+1} = \Pi_{\alpha \leq \beta} X_\alpha = X_\beta \times X^\beta$. Let $\pi^\beta$ be the projection from $X^{\beta+1}$ along $X_\beta$ onto $X^\beta$ and for $S$ a subset of $X$ set $S^\beta = \pi^\beta(S)$.

**1.2 Theorem** (Noble [1971, Theorem 1.8]). Let $X = \Pi_{\alpha \in \kappa} X_\alpha$ where $\kappa$ is a limit ordinal. If for each $\beta$ satisfying $0 < \beta < \kappa$ the projection of $X^{\beta+1}$ onto $X^\beta$ is closed, then for each $\beta < \kappa$ the projection of $X$ onto $X^\beta$ is closed.

**Proof.** Suppose $\beta < \kappa$ and $S \subseteq X$ are such that $S^\beta$ is not closed; by setting $X^\beta$ equal to $X_0$ and, for $\alpha > \beta$, re-indexing $X_\alpha$ as $X_{\alpha - \beta}$, we may suppose $\beta = 0$. Choose a point $x_0$ in $\text{cl}(S^0) \setminus S^0$ and for $\beta > 0$ exploit Proposition 1.1 to choose inductively points $x_\beta$ in $X_\beta$ such that $(x_0,...,x_\beta)$ is in $\text{cl}(S^{\beta+1})$. Notice that for each $\beta \leq \kappa$ the point $(x_\alpha : \alpha < \beta)$ is in the closure of $S^\beta$, either by construction or, for $\beta$ a limit ordinal, because verification that a basic open neighborhood of a point meets a subset of a product depends on only finitely many coordinates. In particular, the point $(x_\alpha : \alpha < \kappa)$ is in $\text{cl}(S)$. Since $x_0$ is not in $S^0$ this point is not in $S$, $S$ is not closed, and therefore the projection of $X$ onto each $X^\beta$ is closed.

**1.3 Corollary** (Tyconoff's Theorem). Each product of compact spaces is compact.
**Proof.** Let $X$ be the product of a collection $\{X_\alpha : 0 < \alpha < \kappa\}$ of compact spaces with $X_0$ any space and $\kappa$ any non-zero cardinal, and note that since the projection along each $X_\alpha$ is closed, Theorem 1.2 implies that the projection along $X$ onto $X_0$ is closed. Hence by Mrówka's Theorem $X$ is compact.

## 2 The second proof

Recently I used Theorem 1.2 to prove a more powerful result, repeated below.

**2.1 Theorem** (Noble [2019, Theorem 6.2]).

Let $\kappa$ be an infinite cardinal, suppose **R** is a class of topological spaces closed under products of fewer than $\kappa$ factors, and define **P** and **Q** as follows:
- **P** $= \{X : \pi_Y : X \times Y \to Y$ is closed for each $Y$ in **R**$\}$; and
- **Q** $=$ **P** $\cap$ **R**.

If $X = \Pi_{\alpha \in A} X_\alpha$ where each $X_\alpha$ is in **Q**, then
(a) if $|A| < \kappa$, $X$ is in **Q**; and
(b) if $|A| = \kappa$, $X$ is in **P**.

**Proof.** First note that the class **Q** is closed under finite products; indeed, for $X \in$ **P**, and $Y, Z \in$ **R**, $\pi_{Y \times Z} : X \times Y \times Z \to Y \times Z$ is closed, and for $Y \in$ **P**, $Z \in$ **R**, $\pi_Z : Y \times Z \to Z$ and thus $\pi : X \times Y \times Z \to Z$ is closed. Since this is true for each $Z \in$ **R**, $X \times Y$ is in **P**. Since $X, Y \in$ **Q** implies also that $X \times Y$ is in **R**, $X \times Y$ is in **Q**. To complete the proof, note that $X_0 = X^{\leq 0}$ is in **Q** and suppose inductively that $X^{< \beta}$ is in **Q** for all $\beta < \alpha$. If $\alpha = \beta + 1$ for some such $\beta$, then $X^{\leq \alpha}$ is in **Q** as a finite product of members. Otherwise each $\pi : X^{\leq \alpha} \to X^{< \alpha}$ is closed since $X_\alpha \in$ **P** and $X^{< \alpha} \in$ **R**, so, by Theorem 2.1, $X$ is in **P**. If $|\alpha| < \kappa$ then, as a product of fewer than $\kappa$ factors, $X$ is also in **R** and thus **Q**.

Theorem 2.1 provides a process, given a judicious choice of the class **R**, for generating product theorems involving $\kappa$-fold products: Vaughan [1975, Theorem 2] provides several examples. The limits of that process are suggested by two standard arguments reviewed below.

**2.2 Remark** For spaces X and Y let $\pi_Y: X \times Y \to Y$ be the projection and let $\kappa$ be an infinite cardinal.

(a) If X has a $\kappa$-fold open cover $\{U_\alpha : \alpha \in \kappa\}$ with no subcover of smaller cardinality and Y contains an increasing collection $\{H_\alpha : \alpha \in \kappa\}$ of closed sets whose union is not closed, then $\pi_Y$ is not closed.

(b) If X is $[\kappa, \infty]$-compact (each open cover has a subcover of cardinality less than $\kappa$) and in Y each intersection of fewer than $\kappa$ open sets is open, then $\pi_Y$ is closed.

(c) If $\pi_Y$ is closed and Y is a product with $\kappa$ nontrivial (not indiscrete) factors, then X is initially $\kappa$-compact (each $\kappa$-fold or smaller open cover has a finite subcover).

**Proof.** (a). Set $F_\alpha = \cap_{\beta < \alpha}(X \setminus U_\beta)$; then $S = \cup_{\alpha < \kappa} F_\alpha \times H_\alpha$ is closed but $\pi_Y S = \cup_\alpha H_\alpha$ is not.

(b). If S is a closed subset of X x Y and y is a point of Y not in $\pi_Y S$, we can cover X x {y} with $\kappa$-many basic open subsets $\{U_\alpha \times V_\alpha : \alpha \in \kappa\}$ which do not meet S, demonstrating that $\cap_\alpha V_\alpha$ is a neighborhood of y which does not meet $\pi_Y S$. Hence $\pi_Y S$ is closed and therefore $\pi_Y$ is closed.

(c). Suppose $Y = \Pi_{\alpha \in \kappa} Y_\alpha$ where each $Y_\alpha$ contains a non-empty proper closed subset $C_\alpha$ and set $H_\alpha = \Pi_{\beta \leq \alpha} X_\beta \times \Pi_{\alpha < \beta} C_\beta$. Note that if $\lambda \leq \kappa$ is an infinite cardinal, $\{H_\alpha : \alpha \in \lambda\}$ is an increasing collection of closed sets whose union is not closed; thus by (a) and the fact that $\pi_Y$ is closed any infinite cardinal $\lambda$ such that X has a $\lambda$-fold open cover with no subcover of smaller cardinality must be greater than $\kappa$. It follows that X is initially $\kappa$-compact.

In particular, if X is a Lindelöf space which is not compact then for any space Y the projection $\pi_Y$ on X x Y is closed if and only if Y is a P-space. Notice that a product of Lindelöf P-spaces with a countably infinite number of nontrivial factors is not a P-space, so 2.1(b) cannot be improved to conclude that X is in **Q**. Of course, our interest here is with the application of 2.1 to the simpler consequence of (b) and (c): Mrówka's Theorem.

**2.3 Corollary** (Tyconoff's Theorem yet again). Each product of compact spaces is compact.

**Proof.** In Theorem 2.1, take **R** to be class of topological spaces.